\documentclass[12pt,a4paper,oneside,reqno]{article}

\usepackage[utf8]{inputenc}
\usepackage[english]{babel}
\usepackage[left=1in,right=1in,top=1in,bottom=1in]{geometry}
\usepackage{amsfonts,amssymb,amsmath,amsthm}
\usepackage[symbol]{footmisc}
\usepackage{longtable,array}
\newcolumntype{L}[1]{>{\raggedright\arraybackslash$}p{#1}<{$}}
\usepackage{enumerate}
\usepackage{cite}
\usepackage{xcolor}

\usepackage{hyperref}
\hypersetup{hidelinks}
\usepackage{fancyhdr}
\theoremstyle{plain}
\newtheorem{theorem}{Theorem}

\theoremstyle{definition}

\makeatletter
\@removefromreset{theorem}{section}\@removefromreset{equation}{section}

\makeatother
\usepackage{adjustbox}
\usepackage{tikz}\usetikzlibrary{arrows}

\DeclareMathOperator{\tr}{tr}
\DeclareMathOperator{\Der}{Der}

\IfFileExists{fouriernc.sty}{\usepackage{fouriernc}}{}

\newcommand{\doi}[1]{\href{https://doi.org/#1}{doi:\detokenize{#1}}}
\newcommand{\arxiv}[1]{\href{https://arxiv.org/abs/#1}{arXiv:#1}}

\theoremstyle{plain}
\newtheorem{Th}[theorem]{Theorem}
\newtheorem{Lem}[theorem]{Lemma}
\theoremstyle{definition}

\newtheorem{Example}[theorem]{Example}
\newcommand{\FF}{{\mathbb{F}}}
\newcommand{\KK}{{\mathbb{K}}}
\renewcommand{\tr}{\mathop{\rm Tr}}

\begin{document}

\noindent\fbox{\parbox{\dimexpr\textwidth-2\fboxsep-2\fboxrule\relax}{\small
\textbf{Corrected version (verified by Claude).} This is a corrected version of the paper I.~Kaygorodov, A.~Lopatin, and Yu.~Popov, \emph{Conservative algebras of $2$-dimensional algebras}, Linear Algebra and its Applications, 486 (2015), 255--274, \doi{10.1016/j.laa.2015.08.011}.}}

\bigskip

\begin{center}
{\Large Conservative algebras of $2$-dimensional algebras}

\bigskip

{\bf Ivan Kaygorodov, Artem Lopatin, and Yury Popov}
\end{center}

\medskip

\noindent{\small\textbf{Abstract.} In 1990 Kantor introduced the conservative algebra $W(n)$ of all algebras (i.e. bilinear maps) on the $n$-dimensional vector space. In the case $n>1$, the algebra $W(n)$ does not belong to any of the well-known classes of algebras (such as associative, Lie, Jordan, or Leibniz algebras). We describe the algebra of all derivations of $W(2)$ and the subalgebras of $W(2)$ of codimension one. We also study similar problems for the algebra $W_2$ of all commutative algebras on the two-dimensional vector space and the algebra $S_2$ of all commutative algebras with trace zero multiplication on the two-dimensional space.}

\bigskip

\noindent {\bf Keywords:}
{\it bilinear map, conservative algebra, terminal algebra, derivation, subalgebra of codimension one.}

\bigskip

\noindent {\bf MSC2020:} 17A30 (primary); 15A03, 15A69, 17A36, 17D10 (secondary).

\bigskip

\tableofcontents

\bigskip

\section*{Introduction}\addcontentsline{toc}{section}{Introduction}
We work over an arbitrary field $\FF$ of characteristic zero. An algebra may not have a unity and may not be associative. By a multiplication on a vector space $W$ we mean a bilinear mapping $W\times W \to W$. Given a vector space $W$, a linear mapping $A$ on $W$, and a bilinear mapping $B:W\times W\to W$, we can define a multiplication $[A,B]:W\times W\to W$ as follows:
$$[A,B](x,y)=A(B(x,y))-B(A(x),y)-B(x,A(y))$$
for all $x,y\in W$. For an algebra $A$ with multiplication $P$ and $x\in A$ we write $L_x^P$ for the operator of left multiplication by $x$. If the multiplication $P$ is fixed, then we write $L_x$ instead of $L_x^P$.

In 1990 Kantor~\cite{Kantor90} defined the multiplication $\cdot$ on the set of all algebras (i.e. all multiplications) on the $n$-dimensional vector space $V_n$ as follows:
$$A\cdot B = [L_{e}^A,B],$$
where $A$ and $B$ are multiplications and $e$ is the fixed vector from $V_n$. In the case $n>1$, the resulting algebra $W(n)$ does not belong to any of the well-known classes of algebras (such as associative, Lie, Jordan, Leibniz algebras). The algebra $W(n)$ turns out to be a conservative algebra (see below).

In 1972 Kantor~\cite{Kantor72} introduced conservative algebras as a generalization of Jordan algebras. Namely, an algebra $A$ with the multiplication $P(x,y)=xy$ and the underlying vector space $W$ is a {\it conservative algebra} if  a new multiplication $F$  can be defined on the vector space $W$ in such a way that
\begin{eqnarray}\label{tojd_oper}
[L_b, [L_a,P]]=-[L_{F(a,b)},P]
\end{eqnarray}
for all $a,b\in W$. In other words, the following identity holds for all $a,b,x,y$ from $W$:
\begin{multline*}
b(a(xy)-(ax)y-x(ay))-a((bx)y)+(a(bx))y\\+(bx)(ay)-a(x(by))+(ax)(by)+x(a(by))=
\end{multline*}
\begin{eqnarray}
-F(a,b)(xy)+(F(a,b)x)y+x(F(a,b)y).\label{tojdestvo_glavnoe}
\end{eqnarray}

\noindent{}The algebra with the multiplication $F$ is said to be {\it associated} to $A$.

As an example, it is easy to see that every $4$-nilpotent algebra is a conservative algebra with $F(a,b)=0$.

The algebra $W(n)$ plays a similar role in the theory of conservative algebras as the Lie algebra of all $n\times n$ matrices $gl_n$ plays in the theory of Lie algebras. Namely, in~\cite{Kantor90}  Kantor considered the category $\mathcal{S}_n$ whose objects were conservative algebras of non-Jacobi dimension $n$. It was proven that the algebra $W(n)$ is the universal attracting object in this category, i.e., for every algebra $M$ of $\mathcal{S}_n$ there exists a canonical homomorphism from $M$ into the algebra $W(n)$. In particular, all Jordan algebras of dimension $n$ with unity are contained in the algebra $W(n)$. The same statement also holds for all noncommutative Jordan algebras of dimension $n$ with unity.
Some properties of the product in the algebra $W(n)$ were studied in \cite{kaykantor}.

Let us consider the known results on conservative algebras. Kantor classified all conservative algebras of order two in~\cite{Kantor72} and defined the class of {\it terminal} algebras as algebras satisfying a certain identity. He proved that every terminal algebra is a conservative algebra and classified all simple finite-dimensional terminal algebras with a left quasi-unit over an algebraically closed field of characteristic zero~\cite{Kantor89term}. Terminal trilinear operations were studied in~\cite{Kantor89tril}, and some questions concerning classification of simple conservative algebras were considered in~\cite{Kantor89trudy}. After that, Cantarini and Kac classified simple finite-dimensional (and linearly compact) super-commutative and super-anticommutative conservative superalgebras and some generalization of these algebras (also known as ``rigid''{} superalgebras) over an algebraically closed field of characteristic zero (see~\cite{kac_can_10}).

In Section 2 we summarize the known results about conservative algebras, and construct new examples of conservative and simple non-conservative algebras.

Recently a great interest has been shown to the study of Jordan and Lie algebras and superalgebras, as well as their generalizations with derivations.
Namely, Popov determined the structure of differentiably simple Jordan algebras \cite{Popov2};
Barreiro, Elduque, and Mart\'inez described the derivations of the Cheng--Kac Jordan superalgebra~\cite{BEC};
Kaygorodov and Okhapkina found all $\delta$-derivations of semisimple structurable algebras~\cite{kay_okh}; Kaygorodov, Shestakov, and Zusmanovich studied generalized derivations of Jordan and Lie algebras and superalgebras in~\cite{kay, kay_de2, kay_lie2, zus10, kay_de5, shestakov2}. Another important problem is investigation of subalgebras for associative and non-associative algebras. For example, subalgebras of codimension one were studied by Dzhumadildaev~\cite{dzhuma85}, Wilansky~\cite{wilansky75} and others.

In Section~3 we describe the algebra of all derivations of $W(2)$ and subalgebras of $W(2)$ of codimension one. In Sections~4 and~5 we study similar problems for the algebra $W_2$ of all commutative algebras on $V_2$ and for the algebra $S_2$ of all commutative algebras with trace zero multiplication on $V_2$.

Note that all one-sided ideals, all automorphisms, and all idempotents of the algebras $W(2)$, $W_2$, and $S_2$ were described in~\cite{kayvolk}.

\section{Conservative algebras}

\subsection{Associative algebras}
The variety of associative algebras is defined by the identity
$$(ab)c=a(bc).$$
Every associative algebra is conservative (see~\cite{Kantor72}), since we can set $F(a,b)=ab$.

\subsection{Jordan algebras}
The variety of Jordan algebras is defined by the identities $$ab=ba, (a^2b)a=a^2(ba).$$
Every Jordan algebra is conservative (see~\cite{kac77}) with $F(a,b)=ab$.

\subsection{Structurable algebras}
In 1978 Allison~\cite{Alison78} introduced structurable algebras as one of generalizations of unitary Jordan algebras. We can define a structurable algebra $A$ as a unitary algebra with involution $\tau$, satisfying the identity
$$[V_{x,1}, V_{z,w}] = V_{V_{x,1}z,w} - V_{z, V_{1,x}w}, \mbox{ where }
V_{x,y}(z)=(x\tau{(y)})z+(z\tau{(y)})x-(z\tau{(x)})y.$$
If $(A, \tau)$ is a structurable algebra, then a new multiplication $*$ can be defined on $A$ by
$$x*y=xy+y(x-\tau{(x)}).$$
The algebra $A$ with the multiplication $*$ is a conservative algebra~\cite{AH81}.

\subsection{Terminal algebras}
As we have already mentioned, terminal algebras were introduced in~\cite{Kantor72} and it is known that every terminal algebra is conservative (see~\cite{Kantor89term}).

\subsection{Quasi-associative algebras}
An algebra is called quasi-associative (see, for example,~\cite{dedkov89}) if it satisfies the following identities for a fixed element $\alpha\in\FF$ with $4\alpha+1\neq0$:
$$(x,y,z)+(y,z,x)+(z,x,y)=0,$$
$$(x,y,z) = \alpha [y,[x,z]].$$
Here $(x,y,z)=(xy)z-x(yz)$ is the associator and $[x,y]=xy-yx$ is the commutator. It is known \cite{dedkov89} that an algebra $A$ is quasi-associative if and only if there exists an element $\lambda$ from some  extension $\KK$ of the field $\FF$ with $2\lambda - 1\neq0$ such that the $\KK$-algebra $A\otimes_{\FF}\KK$ with the multiplication
$$x\cdot y = \lambda xy + (1-\lambda) yx$$
is associative. Moreover, we can assume that $\alpha = \lambda (\lambda -1)$. On the other hand, in this case for $\mu=\lambda/(2\lambda - 1)$ we obtain that
$$xy = \mu\; x\cdot y + (1-\mu)\; y\cdot x.$$

For a quasi-associative algebra $A$ with a multiplication $P(x,y)=xy$ it is easy to see that
$[L_a,P](x,y)=-\mu\; x\cdot a\cdot y-(1-\mu)\; y\cdot a\cdot x$.
Hence,
$$[L_b, [L_a,P]](x,y)=
\mu (1-\mu)(x\cdot b\cdot a\cdot y+y\cdot a\cdot b\cdot x)+\mu^2\; x\cdot a\cdot b\cdot y+(1-\mu)^2 y\cdot b\cdot a\cdot x=$$
$$\mu\; x\cdot (\mu\; a\cdot b + (1-\mu)b\cdot a)\cdot y+(1-\mu)\cdot y\cdot (\mu\; a\cdot b+(1-\mu)b\cdot a)\cdot x=$$
$$\mu\; x \cdot(ab)\cdot y +(1-\mu) y \cdot(ab)\cdot x=-[L_{ab}, P](x,y).$$
Assuming that $F(a,b)=ab$ in~(\ref{tojd_oper}) we complete the proof of the following lemma.

\medskip
\begin{Lem} Every quasi-associative algebra is conservative.
\end{Lem}

\medskip

\subsection{Lie algebras}
The variety of Lie algebras is defined by the identities
$$ab=-ba, (ab)c+(bc)a+(ca)b=0.$$
Every Lie algebra is conservative (see~\cite{Kantor72}), since we can set $F(a,b)=0$.

\subsection{Left Leibniz algebras}
The variety of left Leibniz algebras (see, for example, \cite{Dzhuma01})
is defined by the identity $$a(bc)=(ab)c+b(ac).$$
Every left Leibniz algebra is conservative, since $[L_a,P]=0$, where $P$ is the multiplication in this algebra.

\subsection{Malcev algebras}
Malcev algebras were introduced by Malcev in 1955~\cite{Malcev55} as algebras satisfying the identities:
$$xy=-yx,
J(a,x,y)a=J(a,x,ay), \mbox{ where } J(a,x,y)=(ax)y+(xy)a+(ya)x.$$

The best-known example of non-Lie Malcev algebra is a simple $7$-dimensional
Malcev algebra $M_7$. In the case of an algebraically closed field of characteristic zero we can take the following set as a basis for the algebra $M_7$:
$$B_{M_7}=\{ h,x,y,z,x',y',z'\}$$
with the multiplication table
$$hx=2x, hy=2y, hz=2z,$$
$$hx'=-2x',hy'=-2y',hz'=-2z',$$
$$xx'=yy'=zz'=h,$$
$$xy=2z',yz=2x',zx=2y',$$
$$x'y'=-2z,y'z'=-2x,z'x'=-2y,$$
where any other product of elements from the basis is zero  (see~\cite{kuzmin89}).

\medskip
\begin{Lem}\label{lemma_M7}
Over an algebraically closed field the simple $7$-dimensional Malcev algebra $M_7$ is not a conservative algebra.
\end{Lem}

\noindent{\bf Proof.}
Assume that the algebra $M_7$ is conservative. Then it satisfies
identity~(\ref{tojdestvo_glavnoe}). Since $M_7$ is anticommutative, the right-hand side of~(\ref{tojdestvo_glavnoe}) equals $J(F(a,b),x,y)$. We set $a=b=x'$, $x=y'$, $y=z'$ in~(\ref{tojdestvo_glavnoe}). Since $J(x',y',z')=-6h$, a direct computation shows that the left-hand side of~(\ref{tojdestvo_glavnoe}) equals $3J(x',y',z')x'=36x'$. Therefore, $$J(F(a,a),y',z')=36x' \quad \text{for } a=x'.$$
On the other hand,
$$(F(a,a) y') z' \in \langle x,z', h \rangle, \quad F(a,a)(z'y') \in \langle
x,h,z',y' \rangle,\quad \text{and}\; (F(a,a)z')y' \in \langle x,y',h \rangle,$$
where $\langle u,\ldots,v\rangle$ stands for the linear span of vectors $u,\ldots,v$. Hence the projection of the vector $J(F(a,a), y',z')$ on the subspace generated
by $x'$ is equal to zero. It follows that $J(F(a,a),y',z')\neq 36x'$; a contradiction. Thus $M_7$ is not a conservative algebra and the lemma is proven.

\subsection{Non-commutative Jordan algebras}
The class of non-commutative Jordan algebras is defined by the identities
$$(a^2b)a=a^2(ba), (ab)a=a(ba).$$
It includes
all associative algebras, alternative algebras, Jordan algebras, quasi-associative algebras and
all anticommutative algebras (in particular, Lie, Malcev, binary Lie algebras). Every conservative algebra with a unit is a non-commutative Jordan algebra~\cite{Kantor72}. Note that the mentioned identity $(ab)a=a(ba)$ is called the flexibility identity and any algebra, where it holds, is called flexible. It was also shown in~\cite{Kantor72} that every flexible conservative algebra with $F(a,b)=ab$ is a non-commutative Jordan algebra.
An example of a simple non-conservative non-commutative Jordan algebra is the Malcev algebra $M_7$.

\subsection{Poisson algebras}
Let $A$ be a vector space with an associative commutative multiplication $ab$ and a Lie multiplication $\{a,b\}$.
Then the algebra $A$ is called a Poisson algebra if it satisfies the identity
$$ \{ab,c\} =a\{b,c\}+\{a,c\}b.$$

\medskip
\begin{Lem}\label{lem_poisson}
Every Poisson algebra $A$ with a new multiplication $P(a,b)=ab+\{a,b\}$ is a conservative algebra.
\end{Lem}
\noindent{\bf Proof.}
It is known that the algebra $A$ with the multiplication $P$ is a non-commutative Jordan algebra.
Note that
$$[L_a,P](x,y)=-a(xy+\{x,y\})-\{x,a\}y-x\{a,y\},$$
where $L_a(x)=P(a,x)$. Hence,
$$[L_b,[L_a,P]](x,y)=
-(baxy+ba \{x,y\} +b\{x,a\}y+bx\{a,y\}+\{b,axy\}+$$
$$\{b,a \{x,y\}\} +\{b,\{x,a\}y\}+\{b,x\{a,y\}\}-$$
$$ abxy-a\{bx,y\}-\{bx,a\}y-bx\{a,y\} - a\{b,x\}y-a\{\{b,x\},y\}-\{\{b,x\},a\}y-\{b,x\}\{a,y\} -$$
$$ axby-a\{x,by\}-\{x,a\}by-x\{a,by\}  - ax\{b,y\}-a\{x,\{b,y\}\}-\{x,a\}\{b,y\}-x\{a,\{b,y\}\}) =$$
$$ ab(xy+\{x,y\}) +x\{ab,y\} +y\{x,ab\} +\{a,b\}(xy+\{x,y\})+\{x,\{a,b\}\}y+x\{\{a,b\},y\}=$$
$$- [L_{P(a,b)},P](x,y).$$
Taking $F(a,b)=P(a,b)$ we complete the proof of the lemma.
\medskip

\subsection{Left-commutative algebras}
The variety of left-commutative algebras is defined by the identity $$a(bx)=b(ax).$$
An algebra is called a generalized associative algebra if we can introduce a new multiplication $F$ on the underlying vector space of the algebra in such a way that the following identity holds:
$$a(bx)=F(a,b)x.$$

\begin{Th}\label{theorem_left_com}
Let $A$ be a left-commutative conservative algebra with trivial left annihilator $\{z\in A \mid zA=0\}$.
Then $A$ is a generalized associative algebra and the algebra associated to $A$ is commutative.
\end{Th}

\noindent{\bf Proof.}
Since every conservative left-commutative algebra satisfies identity~(\ref{tojdestvo_glavnoe}), then

\begingroup\footnotesize$$b(a(xy)-(ax)y-x(ay))-a((bx)y)+(a(bx))y+(bx)(ay)-a(x(by))+(ax)(by)+x(a(by))=$$

$$[b(a(xy))-b(x(ay))]+[(bx)(ay)-a((bx)y)]+[(ax)(by)-b((ax)y)]+[x(a(by))-a(x(by))]+(a(bx))y=$$
$$(a(bx))y$$\endgroup{}
and
\begin{eqnarray*}
[x(F(a,b)y)-F(a,b)(xy)]+(F(a,b)x)y=(F(a,b)x)y.
\end{eqnarray*}
Thus the identity
\begin{eqnarray}\label{cons_leftcom}
(a(bx))y=(F(a,b)x)y
\end{eqnarray}
holds in every conservative left-commutative algebra.

Equality~(\ref{cons_leftcom}) implies that the element $a(bx)-F(a,b)x$ lies in the left annihilator of $A$. Hence the algebra $A$ satisfies the identity $F(a,b)x=a(bx)$. Since $F(a,b)x=a(bx)=b(ax)=F(b,a)x$, then the algebra associated to $A$ is commutative. The theorem is proven.

Note that the triviality of the left annihilator cannot be replaced by the triviality of the (two-sided) annihilator. Indeed, the $2$-dimensional algebra $A$ with a basis $\{e,n\}$ and the non-zero products $ee=e+n$, $en=n$ is left-commutative, has trivial annihilator, and is conservative with the associated multiplication given by $F(e,e)=e$, $F(e,n)=F(n,e)=F(n,n)=0$. However, $A$ is not a generalized associative algebra, since $e(ee)=e+2n$, whereas $ue\in\langle e+n\rangle$ for every $u\in A$.

\medskip
\begin{Example}\label{ex_leftcom}
For each $n>1$ there exists a simple left-commutative algebra of dimension $n$, which is not conservative.
\end{Example}

\noindent{\bf Proof.}
Given $n>1$, we consider the simple left-commutative algebra $A$ with a basis $\{e_1, \ldots, e_n\}$ and the multiplication table $e_i \cdot e_j=j e_j$ for all $i,j$. This algebra has earlier appeared in~\cite{shestakov_kuzm}.
Assume that $A$ is conservative. Then identity~(\ref{cons_leftcom}) holds in $A$ (see the proof of Theorem~\ref{theorem_left_com}; note that Theorem~\ref{theorem_left_com} itself cannot be applied, since the left annihilator of $A$ is non-trivial). Hence
$$(e_1\cdot (e_1 \cdot e_1))\cdot e_1=(F(e_1,e_1) \cdot e_1)\cdot e_1.$$
We set $F(e_1,e_1)=\sum \alpha_i e_i$ for some $\alpha_1,\ldots,\alpha_n$ from $\FF$. Since $e_1\cdot(e_1\cdot e_1)=e_1$ and $(\sum \alpha_i e_i)\cdot e_1=(\sum\alpha_i)e_1$, we obtain $\sum \alpha_i=1$. On the other hand,
$$(e_1\cdot (e_1 \cdot e_2))\cdot e_1=(F(e_1,e_1)\cdot e_2)\cdot e_1,$$ that is,
$4e_1=2 (\sum \alpha_i) e_1$. Thus $\sum \alpha_i=2$; a contradiction.
Therefore, the algebra $A$ is not conservative.

\subsection{Zinbiel algebras}
The variety of Zinbiel algebras (see, for example, \cite{Dzhuma07}) is defined by the identity
$$a(bx)=(ab+ba)x.$$
Since every Zinbiel algebra is a left-commutative algebra and satisfies identity~(\ref{cons_leftcom}) for $F(a,b)=ab+ba$, then every Zinbiel algebra is a conservative algebra.

\section{Conservative algebra of all $2$-dimensional algebras}

A multiplication of a $2$-dimensional algebra is defined by a $2\times 2\times 2$ array of structure constants. The classification of $2$-dimensional algebras was given in~\cite{GR11,GR2}.
Consider the space $W(2)$ of all multiplications on the $2$-dimensional space $V_2$ with a basis $v_1,v_2$.  The definition of the multiplication $\cdot$  for the algebra $W(2)$ can be found in  Section~1. Namely, we fix a vector $v_1 \in V_2$ and obtain that
$$(A \cdot B)(x,y)=A(v_1,B(x,y))-B(A(v_1,x),y)-B(x,A(v_1,y)),$$
where $x,y \in V_2$ and $A,B \in W(2)$. The algebra $W(2)$ is conservative (see~\cite{Kantor90}) and the multiplication $F$ on the algebra associated to $W(2)$ can be given by the equality
$$F(A,B)=(1/3)(A^* \cdot B+ \tilde{B} \cdot A),\mbox{ where }A^*=A+A^T, \tilde{B}=2B^T-B, A^T(x,y)=A(y,x).$$
Let
$A(v_i,v_j)= a_{ij}^1 v_1+a_{ij}^2 v_2$ and
$B(v_i,v_j)= b_{ij}^1 v_1+b_{ij}^2 v_2$ for some $a_{ij}^1$, $a_{ij}^2$, $b_{ij}^1$, $b_{ij}^2$ from the field $\FF$. Then
\begin{eqnarray}\label{form_AB}
(A \cdot B)(v_i, v_j)=A(v_1, B(v_i, v_j))-B(A(v_1, v_i), v_j)- B(v_i, A(v_1, v_j))=
\end{eqnarray}
$$A(v_1, b_{ij}^1 v_1+b_{ij}^2 v_2)-B(a_{1i}^1 v_1+a_{1i}^2 v_2, v_j)-B(v_i, a_{1j}^1 v_1+a_{1j}^2 v_2)=$$
$$(a_{11}^1 b_{ij}^1+a_{12}^1 b_{ij}^2 -  a_{1i}^1 b_{1j}^1-a_{1i}^2 b_{2j}^1-a_{1j}^1 b_{i1}^1- a_{1j}^2 b_{i2}^1) v_1+$$
$$(a_{11}^2 b_{ij}^1+ a_{12}^2 b_{ij}^2 - a_{1i}^1 b_{1j}^2- a_{1i}^2 b_{2j}^2- a_{1j}^1 b_{i1}^2- a_{1j}^2 b_{i2}^2) v_2.$$

We consider the multiplications $\alpha_{i,j}^k$ ($i,j,k=1,2$) on $V_2$ defined by the following formula: $\alpha_{i,j}^k(v_t,v_l)=\delta_{it}\delta_{jl} v_k$ for all $t,l$. It is easy to see that $\{ \alpha_{i,j}^k | i,j,k=1,2 \}$ is a basis of the algebra $W(2)$. Applying~(\ref{form_AB}), we can describe the multiplication table of the conservative algebra $W(2)$:

\begingroup\small\setlength{\arraycolsep}{0.5pt}\begin{eqnarray*}
\begin{array}{l l l l}
\alpha_{11}^1 \cdot \alpha_{11}^1 =- \alpha_{11}^1 	 &
\alpha_{12}^1 \cdot \alpha_{11}^1 =- \alpha_{12}^1-\alpha_{21}^1 	 &
\alpha_{11}^2 \cdot \alpha_{11}^1 = \alpha_{11}^2 	 &
\alpha_{12}^2 \cdot \alpha_{11}^1 =0 	    \\

\alpha_{11}^1 \cdot \alpha_{12}^1 = 0	&
\alpha_{12}^1 \cdot \alpha_{12}^1 = - \alpha_{22}^1 	 &
\alpha_{11}^2 \cdot \alpha_{12}^1 = -\alpha_{11}^1+\alpha_{12}^2	&
\alpha_{12}^2 \cdot \alpha_{12}^1 = -\alpha_{12}^1 	    \\

\alpha_{11}^1 \cdot \alpha_{21}^1 = 0 &
\alpha_{12}^1 \cdot \alpha_{21}^1 =- \alpha_{22}^1 	 &
\alpha_{11}^2 \cdot \alpha_{21}^1 =-\alpha_{11}^1+\alpha_{21}^2	 &
\alpha_{12}^2 \cdot \alpha_{21}^1 =-\alpha_{21}^1  	    \\

\alpha_{11}^1 \cdot \alpha_{22}^1 = \alpha_{22}^1 &
\alpha_{12}^1 \cdot \alpha_{22}^1 = 0  	 &
\alpha_{11}^2 \cdot \alpha_{22}^1 =-\alpha_{12}^1-\alpha_{21}^1+\alpha_{22}^2	 &
\alpha_{12}^2 \cdot \alpha_{22}^1 =-2\alpha_{22}^1  	    \\

\alpha_{11}^1 \cdot \alpha_{11}^2 = -2\alpha_{11}^2 &
\alpha_{12}^1 \cdot \alpha_{11}^2 = \alpha_{11}^1-\alpha_{21}^2-\alpha_{12}^2  	 &
\alpha_{11}^2 \cdot \alpha_{11}^2 = 0	 &
\alpha_{12}^2 \cdot \alpha_{11}^2 = \alpha_{11}^2  	    \\

\alpha_{11}^1 \cdot \alpha_{12}^2 = -\alpha_{12}^2 &
\alpha_{12}^1 \cdot \alpha_{12}^2 =  \alpha_{12}^1-\alpha_{22}^2  	 &
\alpha_{11}^2 \cdot \alpha_{12}^2 = -\alpha_{11}^2	 &
\alpha_{12}^2 \cdot \alpha_{12}^2 = 0  	    \\

\alpha_{11}^1 \cdot \alpha_{21}^2 = -\alpha_{21}^2 &
\alpha_{12}^1 \cdot \alpha_{21}^2 =  \alpha_{21}^1-\alpha_{22}^2  	 &
\alpha_{11}^2 \cdot \alpha_{21}^2 = -\alpha_{11}^2	 &
\alpha_{12}^2 \cdot \alpha_{21}^2 = 0  	    \\

\alpha_{11}^1 \cdot \alpha_{22}^2 = 0  &
\alpha_{12}^1 \cdot \alpha_{22}^2 =  \alpha_{22}^1   	 &
\alpha_{11}^2 \cdot \alpha_{22}^2 = -\alpha_{12}^2-\alpha_{21}^2	 &
\alpha_{12}^2 \cdot \alpha_{22}^2 = -\alpha_{22}^2
\end{array}
\end{eqnarray*}\endgroup{}
\medskip

\noindent{}The rest of the products are zero. It is convenient to denote elements of the basis as follows:
$$e_1=\alpha_{11}^1,\, e_2=\alpha_{12}^1,\, e_3=\alpha_{21}^1,\,   e_4=\alpha_{22}^1,\, e_5=\alpha_{11}^2,\, e_6=\alpha_{12}^2,\,  e_7=\alpha_{21}^2,\, e_8=\alpha_{22}^2.$$

In what follows we write $M$ for an algebra. It is well-known that a linear map $D:M\to M$ is called a {\it derivation} if $D(ab)=D(a)b + a D(b)$ for all $a,b$ from $M$. The set of all derivations of $M$ forms the Lie algebra $\Der(M)$ with respect to the multiplication $[D_1,D_2]= D_1 D_2 - D_2 D_1$.

\medskip
\noindent{\bf Definition.}
An element $a$ of an algebra $M$ is called a {\it Jacobi} element if
$$a(xy)=(ax)y+x(ay).$$
In other words, the transformation $L_a$ is a derivation if and only if $a$ is a Jacobi element. Note that all Jacobi elements of $M$ form a vector space, which is  called the Jacobi space of $M$.
\medskip

It follows from~\cite{Kantor90} that the codimension of the Jacobi space $J$ of the algebra $W(2)$ is two. Using the multiplication table of $W(2)$,
we can see that $J$ is the $\FF$-span of $\alpha_{ij}^k$ for all $i+j>2$.

\medskip
\noindent{\bf Definition.}
An element $e$ of $M$ is said to be a {\it left quasiunit} if
$$e(xy)=(ex)y+x(ey)-xy \text{ for all } x,y\in M$$

It is obvious that if $e$ is a left quasiunit, then for any Jacobi element $x$  we have another left quasiunit $e+x$. Kantor~\cite{Kantor90} noted that the algebra $W(2)$ has a left quasiunit.
We look for a left quasiunit $e$ as $\alpha \alpha_{11}^1+\beta \alpha_{11}^2$, where $\alpha,\beta\in \FF$. Straightforward computations provide the equality $e=-\alpha_{11}^1$. It is easy to see that $-\alpha_{11}^1$ is a left  quasiunit, but it is not a  left unit.

\medskip
\begin{Th}\label{thm_derW}
The algebra of derivations of $W(2)$ is the solvable $2$-dimensional Lie algebra with non-zero multiplication.
\end{Th}

\noindent{\bf Proof.}
Given a derivation $D$ of the algebra $W(2)$, we have $D(e_i)= \sum x_{ij} e_j$ for some $x_{ij}\in\FF$. We consider some relations between images of $e_1,\ldots,e_8$ with respect to $D$ to obtain a system of linear equations on $\{x_{ij}\}$.

Since $-D(e_1)=D(e_1^2)=D(e_1)e_1+e_1 D(e_1),$ then
$$x_{11}=0, x_{12}=x_{13}, x_{14}=0, x_{18}=0.$$

Since $0=D(e_1e_2)=D(e_1)e_2+e_1D(e_2),$
we have
$$x_{21}=-x_{15}, x_{25}=0, x_{16}=0, x_{12}=x_{24}, x_{26}=x_{15}, x_{27}=0.$$

Since $0=D(e_1e_3)=D(e_1)e_3+e_1D(e_3),$
we have
$$x_{31}=-x_{15}, x_{34}=x_{12}, x_{35}=0, x_{36}=0, x_{37}=x_{15}.$$

Since $-D(e_2+e_3)=D(e_2e_1)=D(e_2)e_1+e_2D(e_1),$
we have
$$x_{31}=-x_{15}, x_{32}=0, x_{17}=x_{28}+x_{38}.$$

Since $-D(e_4)=D(e_2e_2)=D(e_2)e_2+e_2D(e_2),$
we have
$$x_{41}=0, x_{43}=x_{21}, x_{44}=2x_{22}+x_{23}-x_{28},x_{48}=x_{26}.$$

Since $-D(e_4)=D(e_2e_3)=D(e_2)e_3+e_2D(e_3),$
we have
$$x_{44}=x_{22}+x_{33}-x_{38}.$$

Since $D(e_5)=D(e_5)e_1+e_5D(e_1),$
we have
$$x_{51}=-x_{12}, x_{52}=0, x_{53}=0, x_{54}=0, x_{56}=x_{12}, x_{57}=x_{13}, x_{58}=0, x_{17}=0.$$

Since $-D(e_6)=D(e_1)e_6+e_1D(e_6),$
we have
$$-x_{62}=x_{12}, x_{63}=0, x_{64}=0, x_{65}=-x_{15}, x_{68}=x_{12}.$$

Since $-D(e_7)=D(e_1)e_7+e_1D(e_7),$
we have
$$x_{72}=0, x_{73}=-x_{12}, x_{74}=0, x_{75}=-x_{15}, x_{78}=x_{12}.$$

Since $0=D(e_1)e_8+e_1D(e_8),$
we have
$$x_{81}=0, x_{12}=-x_{84}, x_{85}=0, x_{15}=-x_{86}, x_{87}=-x_{15}.$$

Since $D(e_1-e_6-e_7)=D(e_2e_5)=D(e_2)e_5+e_2D(e_5),$
we have
$$x_{71}+x_{61}=-x_{22}-x_{55}, x_{76}+x_{66}=x_{22}+x_{55}=x_{77}+x_{67}.$$

Since $D(-e_1+e_6)=D(e_5)e_2+e_5D(e_2),$
we have
$$x_{61}=-x_{55}-x_{22}-x_{23}, x_{66}=x_{55}+x_{22}-x_{28}, x_{67}=x_{23}-x_{28} \mbox{ and } x_{28}=x_{76}.$$

Since $D(-e_1+e_7)=D(e_5e_3)=D(e_5)e_3+e_5D(e_3),$
we have
$$x_{71}=-x_{55}-x_{33}, x_{76}=-x_{38}, x_{77}=x_{55}+x_{33}-x_{38}.$$

The equality $D(-e_2-e_3+e_8)=D(e_5e_4)=D(e_5)e_4+e_5D(e_4)$
implies that
$$x_{22}-x_{82}=x_{23}+x_{33}-x_{83}=x_{55}+x_{44}=-x_{28}-x_{38}+x_{88}.$$

The equality $D(e_2 -e_8)=D(e_2e_6)=D(e_2)e_6+e_2D(e_6)$
implies that
$$x_{82}=x_{61}-x_{66}, x_{23}-x_{83}=x_{67}-x_{61}, x_{88}=x_{22}+x_{28}+x_{66}+x_{67}.$$

The equality $D(e_3-e_8)=D(e_2e_7)=D(e_2)e_7+e_2D(e_7)$
implies that
$$x_{82}=-x_{71}+x_{76}, x_{33}+x_{83}=x_{22}-x_{71}+x_{77}, x_{88}=x_{22}+x_{38}+x_{76}+x_{77}.$$

The equality $D(e_4)=D(e_2e_8)=D(e_2)e_8+e_2D(e_8)$
implies that
$$x_{26}+x_{48}+x_{86}+x_{87}=0.$$

The equality $-D(e_5)=D(e_5e_6)=D(e_5)e_6+e_5D(e_6)$
implies that
$$x_{61}=x_{66}+x_{67}.$$

The equality $-D(e_5)=D(e_5e_7)=D(e_5)e_7+e_5D(e_7)$
implies that
$$x_{71}=x_{76}+x_{77}.$$

The equality $-D(e_6+e_7)=D(e_5e_8)=D(e_5)e_8+e_5D(e_8)$
implies that
$$x_{61}+x_{71}=x_{83}+x_{82}, x_{66}+x_{76}=x_{55}-x_{82}+x_{88}, x_{67}+x_{77}=x_{55}-x_{83}+x_{88}.$$

The equality $0=D(e_6e_1)=D(e_6)e_1+e_6D(e_1)$
implies that
$$x_{61}=0.$$

The equality $0=D(e_7e_1)=D(e_7)e_1+e_7D(e_1)$
implies that
$$x_{71}=0, x_{75}=0.$$

Thus, we can easily see that
$$ x_{12}=x_{13}=x_{24}=x_{34}=-x_{51}=x_{56}=x_{57}=-x_{62}=x_{68}=-x_{73}=x_{78}=-x_{84}=z,$$
$$x_{22}=x_{33}= \frac{x_{44}}{2} =-x_{55}=x_{88}=w$$
for some $z,w\in\FF$. Moreover, any other element $x_{ij}$ is equal to zero.
Then the linear map $D$ is given by the following matrix
\begin{equation}\label{table1}\left(\begin{array}{ccccccccc}
0 & z & z & 0 & 0 & 0 & 0  & 0  \\
0 & w & 0 & z & 0 & 0 & 0  & 0  \\
0 & 0 & w & z & 0 & 0 & 0  & 0  \\
0 & 0 & 0 & 2w & 0 & 0 & 0 & 0  \\
-z & 0 & 0 & 0 & -w & z & z & 0  \\
0 & -z & 0 & 0 & 0 & 0 & 0 & z  \\
0 & 0 & -z & 0 & 0 & 0 & 0 & z  \\
0 & 0 & 0 & -z & 0 & 0 & 0 & w
\end{array} \right)
\end{equation}
and the dimension of $\Der(W(2))$ is two.

Note that the left multiplications $L_{e_2}$ and $L_{e_6}$ are derivations of $W(2)$. Since $[L_{e_2},L_{e_6}]=L_{e_2}$, the derivations $L_{e_2}$ and $L_{e_6}$ span the two-dimensional non-abelian Lie algebra $\Der(W(2))$. The theorem is proven.

\medskip

\begin{Th}\label{thm_codimW}
Let $B$ be a subalgebra of $W(2)$ of codimension one. Then $B$ is the linear span of $e_1,e_2,e_3,e_4,e_6,e_7,e_8$.
\end{Th}

\noindent{\bf Proof.}
It follows from the multiplication table that $e_3$, $e_4$, $e_7$, $e_8$ are left annihilators of $W(2)$, that is, $e_3W(2)=e_4W(2)=e_7W(2)=e_8W(2)=0$, and that the operator $L_{e_1}$ is diagonal in the basis $e_1,\ldots,e_8$:
$$L_{e_1}e_1=-e_1,\quad L_{e_1}e_6=-e_6,\quad L_{e_1}e_7=-e_7,\quad L_{e_1}e_4=e_4,$$
$$L_{e_1}e_2=L_{e_1}e_3=L_{e_1}e_8=0,\quad L_{e_1}e_5=-2e_5.$$
Hence $W(2)=V_{-1}\oplus V_{0}\oplus V_{1}\oplus V_{-2}$, where $V_{-1}=\langle e_1,e_6,e_7\rangle$, $V_0=\langle e_2,e_3,e_8\rangle$, $V_1=\langle e_4\rangle$, and $V_{-2}=\langle e_5\rangle$ are the eigenspaces of $L_{e_1}$.

Let $B$ be a subalgebra of $W(2)$ of codimension one. First, we show that $B$ contains an element $h$ with $L_h=L_{e_1}$. If $e_1\in B$, then we take $h=e_1$. Let $e_1\notin B$. Then $W(2)=B\oplus\langle e_1\rangle$, and for every $u\in W(2)$ there is $\lambda\in\FF$ such that $u+\lambda e_1\in B$. If $e_k+\lambda e_1\in B$ with $\lambda\neq0$ for some $k\in\{3,4,7,8\}$, then we take $h=\lambda^{-1}(e_k+\lambda e_1)$. Otherwise $e_3,e_4,e_7,e_8\in B$. Take $b=e_5+\lambda e_1\in B$. Then
$$be_4=-e_2-e_3+e_8+\lambda e_4\in B,\quad be_8=-e_6-e_7\in B,$$
hence $e_2,e_6\in B$, and $e_2b=e_1-e_6-e_7-\lambda(e_2+e_3)\in B$ gives $e_1\in B$; a contradiction.

Since $L_h=L_{e_1}$ and $hB\subseteq B$, the subspace $B$ is invariant under the diagonalizable operator $L_{e_1}$; hence $B=\bigoplus_{\mu}(B\cap V_\mu)$. Since $B$ has codimension one, there is a unique $\mu$ such that $V_\nu\subseteq B$ for all $\nu\neq\mu$ and $B\cap V_\mu$ has codimension one in $V_\mu$.
\begin{enumerate}
\item[1)] If $\mu=-1$, then $e_2,e_3,e_4,e_5,e_8\in B$, and the equalities $e_5e_2=-e_1+e_6$, $e_5e_3=-e_1+e_7$, $e_5e_8=-e_6-e_7$ give $e_1,e_6,e_7\in B$; a contradiction.

\item[2)] If $\mu=0$, then $e_1,e_4,e_5,e_6,e_7\in B$ and $e_5e_4=-e_2-e_3+e_8\in B$. In particular, there is $x=\alpha e_2+\beta e_3+\gamma e_8\in B$ with $\alpha\neq0$, and then $xe_6=\alpha(e_2-e_8)\in B$ and $xe_7=\alpha(e_3-e_8)\in B$. Together with $-e_2-e_3+e_8\in B$, this gives $e_8\in B$, and then $e_2,e_3\in B$; a contradiction.

\item[3)] If $\mu=1$, then $e_2,e_8\in B$ and $e_2e_8=e_4\in B$; a contradiction.

\item[4)] If $\mu=-2$, then $B=\langle e_1,e_2,e_3,e_4,e_6,e_7,e_8\rangle$. This subspace is a subalgebra, since, by the multiplication table, a product of two basis elements has a non-zero coefficient at $e_5$ only if one of the factors is $e_5$.
\end{enumerate}
The theorem is proven.

\section{Terminal algebra $W_2$ of all $2$-dimensional commutative algebras}

Consider the subalgebra $W_2$ of $W(2)$ that was introduced by Kantor in~\cite{Kantor89term,Kantor90} as the algebra of all bilinear commutative operations on the $2$-dimensional space $V_2$ with a basis $v_1,v_2.$ See above for the definition of the multiplication $\cdot$ on $W_2$. The algebra $W_2$ is known to be  terminal~\cite{Kantor89term} and therefore it is conservative. The multiplication $F$ on an algebra associated to $W_2$ can be given by the equality
$$F(A,B)=(1/3)(2 A \cdot B+ B \cdot A).$$

For short, we denote
$\xi_1=\alpha_{11}^1$, $\xi_2=\alpha_{12}^1+\alpha_{21}^1$, $\xi_3=\alpha_{22}^1$, $\xi_4=\alpha_{11}^2$, $\xi_5=\alpha_{12}^2+\alpha_{21}^2$ and $\xi_6=\alpha_{22}^2$, where $\{\alpha_{ij}^k\}$ are multiplications on $V_2$ introduced in Section 3.
It is easy to see that every multiplication $\xi_i$ is a commutative multiplication and the multiplication $\alpha_{12}^1$ is not commutative.
Assume that the dimension of the algebra $W_2$ is 7.  Since $\xi_2$ and $\xi_5$ are elements of $W_2$, we obtain that $\beta \alpha_{12}^1 + \gamma \alpha_{12}^2$ is also a non-zero element of $W_2$ for some $\beta,\gamma\in\FF$.
Hence
$$ \beta v_1 + \gamma v_2= (\beta \alpha_{12}^1 + \gamma \alpha_{12}^2)(v_1,v_2)= (\beta \alpha_{12}^1 + \gamma \alpha_{12}^2)(v_2,v_1)=0.$$
Thus the dimension of the algebra $W_2$ is 6,
and $\{ \xi_i \}$ is a basis of this algebra.
Using~(\ref{form_AB}), we can describe the multiplication table of the terminal algebra $W_2$:

\begin{eqnarray*}
\begin{array}{l l l l}
\xi_1 \cdot \xi_1  =- \xi_1 	 &
\xi_2 \cdot \xi_1  =- \xi_2 	 &
\xi_4 \cdot \xi_1  =  \xi_4 	 &
\xi_5 \cdot \xi_1  = 0 	 \\

\xi_1 \cdot \xi_2  = 0 	 &
\xi_2 \cdot \xi_2  =- 2 \xi_3 	 &
\xi_4 \cdot \xi_2  = \xi_5-2 \xi_1 	 &
\xi_5 \cdot \xi_2  = -\xi_2 	 \\

\xi_1 \cdot \xi_3  = \xi_3 	 &
\xi_2 \cdot \xi_3  = 0 	 &
\xi_4 \cdot \xi_3  = -\xi_2+\xi_6	 &
\xi_5 \cdot \xi_3  = -2\xi_3 	 \\

\xi_1 \cdot \xi_4  = -2\xi_4 	 &
\xi_2 \cdot \xi_4 =  \xi_1 -\xi_5 	 &
\xi_4 \cdot \xi_4  = 0	 &
\xi_5 \cdot \xi_4  = \xi_4 	 \\

\xi_1 \cdot \xi_5  = -\xi_5 	 &
\xi_2 \cdot \xi_5  =  \xi_2  -2\xi_6 	 &
\xi_4 \cdot \xi_5  = -2\xi_4	 &
\xi_5 \cdot \xi_5  = 0 	 \\

\xi_1 \cdot \xi_6  = 0 	 &
\xi_2 \cdot \xi_6  = \xi_3 	 &
\xi_4 \cdot \xi_6  = -\xi_5	 &
\xi_5 \cdot \xi_6  = -\xi_6. 	 \\

\end{array}
\end{eqnarray*}
\medskip

\noindent{}The rest of products are zero.

\begin{Th}\label{thm_derW2}
The algebra of derivations of $W_2$ is isomorphic to the solvable $2$-dimensional Lie algebra with non-zero multiplication.
\end{Th}

\noindent{\bf Proof.}
Let $D$ be a derivation of $W_2$ and $D(\xi_i)=\sum \alpha_{ij} \xi_j$ for some $\alpha_{ij}$ from $\FF$.

Then the equality $-D(\xi_1)=D(\xi_1 \xi_1)=D(\xi_1)\xi_1+\xi_1D(\xi_1)$
implies that $$\alpha_{11}=\alpha_{13}=\alpha_{16}=0.$$

Since $0=D(\xi_1)\xi_2+\xi_1D(\xi_2)$, we can see that
$$\alpha_{15}=\alpha_{24}=0, \alpha_{21}=-2\alpha_{14}, \alpha_{23}=2\alpha_{12}, \alpha_{14}=\alpha_{25}.$$

Since $-D(\xi_2)= D(\xi_2)\xi_1+\xi_2D(\xi_1)$, we obtain
$$\alpha_{14}=\alpha_{24}=\alpha_{25}=\alpha_{26}=0 \text{ and } \alpha_{23}=2\alpha_{12}.$$

Note that
$D(\xi_3)=-\frac{1}{2}(D(\xi_2)\xi_2+\xi_2D(\xi_2))= 2\alpha_{22} \xi_3.$
The equality $D(\xi_4)=D(\xi_4)\xi_1+\xi_4 D(\xi_1)$ implies that
$$\alpha_{41}=-\alpha_{12}, \alpha_{45}=\alpha_{12}, \alpha_{42}=\alpha_{43}=\alpha_{46}=0.$$

It follows from $0=D(\xi_5)\xi_1+\xi_5D(\xi_1)$
that $$\alpha_{51}=\alpha_{54}=0, \alpha_{52}=-\alpha_{12}.$$

Since $-D(\xi_5)=D(\xi_1)\xi_5+\xi_1 D(\xi_5)$,
then $$\alpha_{53}=0, \alpha_{56}=2\alpha_{12}.$$

Since $D(\xi_1-\xi_5)=D(\xi_2)\xi_4+\xi_2 D(\xi_4)$
then $$\alpha_{22}+\alpha_{44}=\alpha_{55}=0.$$

Since $0=D(\xi_1)\xi_6+\xi_1D(\xi_6)$
then
$$\alpha_{61}=\alpha_{64}=\alpha_{65}=0, \alpha_{63}=-\alpha_{12}.$$

Since $D(\xi_2-2\xi_6)=D(\xi_2)\xi_5+\xi_2D(\xi_5)$
then $$\alpha_{66}=\alpha_{22}.$$

Since $D(\xi_3)=D(\xi_2)\xi_6+\xi_2 D(\xi_6)$
then $$\alpha_{62}=0.$$

Therefore, it is easy to see that
$$\alpha_{12}= \frac{\alpha_{23}}{2}=-\alpha_{41}=\alpha_{45}=-\alpha_{52}=\frac{\alpha_{56}}{2}=-\alpha_{63}=\alpha,$$

$$\alpha_{22}=\frac{\alpha_{33}}{2} =-\alpha_{44}=\alpha_{66}=\beta$$
and any other element $\alpha_{ij}$ is equal to zero. In other words, every element of $\Der(W_2)$ can be given as a linear mapping with the following matrix:

\begin{equation}\label{table2}\left(\begin{array}{ccccccc}
0 & \alpha & 0 & 0 & 0 & 0  \\
0 & \beta & 2\alpha & 0 & 0 & 0  \\
0 & 0 & 2\beta & 0 & 0 & 0  \\
-\alpha  &  0 & 0 & -\beta & \alpha & 0   \\
0 & -\alpha  & 0 & 0 & 0 & 2\alpha  \\
0 & 0 & -\alpha & 0 & 0 & \beta  \\
\end{array} \right).
\end{equation}

Note that the left multiplications $L_{\xi_2}$ and $L_{\xi_5}$ are derivations.
Since  $[L_{\xi_2}, L_{\xi_5}]=L_{\xi_2}$, these two derivations generate a subalgebra of $\Der(W_2)$ of dimension two. The theorem is proven.

\medskip

\begin{Th}\label{thm_codimW2}
Let $B$ be a subalgebra of $W_2$ of codimension $1$.
Then $B$ is the linear span of $\xi_1, \xi_2, \xi_3, \xi_5, \xi_6$.
\end{Th}

\noindent{\bf Proof.}
By the multiplication table, $\xi_3$ and $\xi_6$ are left annihilators of $W_2$, and the operator $L_{\xi_1}$ is diagonal: $L_{\xi_1}\xi_1=-\xi_1$, $L_{\xi_1}\xi_5=-\xi_5$, $L_{\xi_1}\xi_2=L_{\xi_1}\xi_6=0$, $L_{\xi_1}\xi_3=\xi_3$, $L_{\xi_1}\xi_4=-2\xi_4$. Hence $W_2=V_{-1}\oplus V_0\oplus V_1\oplus V_{-2}$ with $V_{-1}=\langle\xi_1,\xi_5\rangle$, $V_0=\langle\xi_2,\xi_6\rangle$, $V_1=\langle\xi_3\rangle$, and $V_{-2}=\langle\xi_4\rangle$.

Let $B$ be a subalgebra of $W_2$ of codimension one. As in the proof of Theorem~\ref{thm_codimW}, we first show that $B$ contains an element $h$ with $L_h=L_{\xi_1}$. If $\xi_1\in B$, then we take $h=\xi_1$. Let $\xi_1\notin B$. If $\xi_k+\lambda\xi_1\in B$ with $\lambda\neq0$ for some $k\in\{3,6\}$, then we take $h=\lambda^{-1}(\xi_k+\lambda\xi_1)$. Otherwise $\xi_3,\xi_6\in B$; take $b=\xi_4+\lambda\xi_1\in B$. Then $b\xi_3=-\xi_2+\xi_6+\lambda\xi_3\in B$ and $b\xi_6=-\xi_5\in B$, so $\xi_2,\xi_5\in B$, and $\xi_2b=\xi_1-\xi_5-\lambda\xi_2\in B$ gives $\xi_1\in B$; a contradiction.

Hence $B$ is invariant under $L_{\xi_1}$, $B=\bigoplus_\mu(B\cap V_\mu)$, and there is a unique $\mu$ such that $V_\nu\subseteq B$ for all $\nu\neq\mu$ and $B\cap V_\mu$ has codimension one in $V_\mu$.
\begin{enumerate}
\item[1)] If $\mu=-1$, then $\xi_2,\xi_3,\xi_4,\xi_6\in B$, and the equalities $\xi_2\xi_4=\xi_1-\xi_5$, $\xi_4\xi_6=-\xi_5$ give $\xi_1,\xi_5\in B$; a contradiction.

\item[2)] If $\mu=0$, then $\xi_1,\xi_3,\xi_4,\xi_5\in B$ and $\xi_4\xi_3=-\xi_2+\xi_6\in B$, so $B\cap V_0=\langle \xi_2-\xi_6\rangle$. However, $(\xi_2-\xi_6)\xi_5=\xi_2-2\xi_6\in B\cap V_0$; a contradiction.

\item[3)] If $\mu=1$, then $\xi_2,\xi_6\in B$ and $\xi_2\xi_6=\xi_3\in B$; a contradiction.

\item[4)] If $\mu=-2$, then $B=\langle \xi_1, \xi_2, \xi_3, \xi_5, \xi_6 \rangle$. This subspace is a subalgebra, since a product of two basis elements has a non-zero coefficient at $\xi_4$ only if one of the factors is $\xi_4$.
\end{enumerate}
The theorem is proven.

\section{Terminal subalgebras $S_2$ and $H_1$ of the algebra $W_2$.}

By definition, the algebra $S_2$ is the subalgebra of $W_2$ consisting of all multiplications $A:V_2\times V_2\to V_2$
such that
$$\tr(L_a^A)=0, \forall \, a\in V_2,$$
where $\tr$ stands for the trace of linear map and $L_a^A(x)=A(a,x)$. Let us recall that $\xi_1,\ldots,\xi_6$ is a basis for $W_2$ and $v_1,v_2$ is a basis for $V_2$. In~\cite{Kantor89term} it was shown that the algebra $S_2$ is terminal. Given a commutative multiplication
$A= x_1 \xi_1+x_2 \xi_2+x_3 \xi_3+ x_4 \xi_4+x_5 \xi_5+x_6 \xi_6$ and an element $a=\alpha_1 v_1 +\alpha_2 v_2$ of $V_2$, where $x_1,\ldots,x_6$, $\alpha_1,\alpha_2$ are elements of $\FF$, we can see that
$\tr(L_a^A)=\alpha_1 x_1+ \alpha_2 x_2+\alpha_1x_5+\alpha_2 x_6$.
Thus $A\in S_2$ if and only if $x_1=-x_5, x_2=-x_6$. Therefore the algebra $S_2$ is the $\FF$-span of $\xi_1-\xi_5$, $\xi_2-\xi_6$, $\xi_3$, $\xi_4$.

For short, we denote
$z_1=\xi_1-\xi_5, z_2=\xi_2-\xi_6, z_3=\xi_3, z_4=\xi_4.$ Obviously, $z_1,\ldots,z_4$ is a basis for $S_2$.
Using~(\ref{form_AB}), we can describe the multiplication table of the terminal algebra $S_2$:

\begin{eqnarray*}
\begin{array}{l l l l}
z_1 \cdot z_1  =- z_1 	 &
z_2 \cdot z_1  =- 2z_2 	 &
z_4 \cdot z_1  =  3z_4 	 \\

z_1 \cdot z_2  = z_2 	 &
z_2 \cdot z_2  =- 3z_3 	 &
z_4 \cdot z_2  =- 2z_1 	 \\

z_1 \cdot z_3  = 3z_3 	 &
z_2 \cdot z_3  = 0 	 &
z_4 \cdot z_3  =- z_2 	 \\

z_1 \cdot z_4  =- 3z_4 	 &
z_2 \cdot z_4  =  z_1 	 &
z_4 \cdot z_4  =  0	 \\
\end{array}
\end{eqnarray*}
\medskip

\noindent{}The rest of products are zero.

By definition, the algebra $H_1$ is the subalgebra of $W_2$ consisting of all multiplications $A$ on $V_2$
preserving the nondegenerate skew-symmetric bilinear form $<,>$ on $V_2$. Namely,
$$<A(x,y),z>+<y,A(x,z)>=0$$
for all $x,y,z$ from $V_2$. Consider a multiplication
$A= x_1 \xi_1+x_2 \xi_2+x_3 \xi_3+ x_4 \xi_4+x_5 \xi_5+x_6 \xi_6$ of $H_1$ for some $x_1,\ldots,x_6\in\FF$. Then
$$<A(v_i,v_i),v_j> = <A(v_i,v_j),v_i> =<A(v_j,v_i),v_i>$$
for all $i,j=1,2$. For $i=1$ and $j=2,$ we obtain that $x_1=-x_5.$
For $i=2$ and $j=1,$ we can see that $x_2=-x_6.$
It is easy to see that multiplications $\xi_1-\xi_5, \xi_2-\xi_6, \xi_3,\xi_4$ preserve the non-degenerate skew-symmetric bilinear form $<,>$, but
the multiplications $\xi_1,\xi_2$ do not preserve it.
Assume that the dimension of the algebra $H_1$ is more than $4$. Then $\alpha \xi_1 +\beta \xi_2$ is a non-zero element of $H_1$ for some $\alpha,\beta\in\FF$.
By straightforward calculations we obtain that $\alpha=\beta=0$; a contradiction. Hence the algebra $H_1$ has the basis $\{z_1,\ldots,z_4\}$ and the algebras $H_1$ and $S_2$ coincide with each other.

\medskip

\begin{Th}\label{yury}
The algebra of derivations of $S_2$ is isomorphic to  the solvable  $2$-dimensional Lie algebra with non-zero multiplication.
\end{Th}

\noindent{\bf Proof.}
Let $D(z_i)=\sum \alpha_{ij}z_j$ be a derivation of $S_2$, where $\alpha_{ij}\in\FF$.

Since $-D(z_1)=D(z_1z_1)=D(z_1)z_1+z_1D(z_1)$, then  $D(z_1)= \alpha_{12}z_2.$
It is easy to see that the left multiplication $L_{z_2}$ is a derivation on $S_2.$
Replacing $D$ by $D+\frac{\alpha_{12}}{2}L_{z_2}$, we may assume that $D(z_1)=0$.
Then $D(z_2)=D(z_1z_2)=z_1D(z_2)$, hence $D(z_2)=\alpha_{22} z_2.$

Since $-3D(z_3)=D(z_2)z_2+z_2D(z_2),$ then $D(z_3)=2\alpha_{22}z_3$.

Since $3D(z_4)=D(z_4)z_1,$ then $D(z_4)=\alpha_{44}z_4.$

Since $0=D(z_1)=D(z_2z_4)=D(z_2)z_4+z_2D(z_4),$ then $\alpha_{22}+\alpha_{44}=0.$

Note that the restriction of the left multiplication $L_{\xi_5}$ (which is a derivation of $W_2$) to the subalgebra $S_2$ is a derivation of $S_2$: it maps $z_1\mapsto 0$, $z_2\mapsto -z_2$, $z_3\mapsto -2z_3$, $z_4\mapsto z_4$. Hence $\Der(S_2)$ is spanned by $L_{z_2}$ and $L_{\xi_5}|_{S_2}$.
Since $[L_{z_2}, L_{\xi_5}|_{S_2}]=L_{z_2}$,
these two derivations span the two-dimensional non-abelian Lie algebra $\Der(S_2)$.
In other words,
every element of $\Der(S_2)$ can be given as a linear mapping with the following matrix (as in~(\ref{table1}) and~(\ref{table2}), the $i$-th row contains the coordinates of the image of $z_i$):

\begin{equation}\label{table3}\left(\begin{array}{ccccccc}
0 & -2\beta & 0 & 0  \\
0 & \alpha & -3\beta & 0  \\
0 & 0 & 2\alpha & 0   \\
\beta  &  0 & 0 & -\alpha    \\
\end{array} \right).
\end{equation}

The theorem is proven.

\medskip

\begin{Th}\label{thm_codimS2}
Let $B$ be a subalgebra of $S_2$ of codimension 1.
Then $B$ is the linear span of $z_1,z_2,z_3$.
\end{Th}

\noindent{\bf Proof.}
By the multiplication table, $z_3$ is a left annihilator of $S_2$, and the operator $L_{z_1}$ is diagonal with pairwise distinct eigenvalues: $L_{z_1}z_1=-z_1$, $L_{z_1}z_2=z_2$, $L_{z_1}z_3=3z_3$, $L_{z_1}z_4=-3z_4$.

Let $B$ be a subalgebra of $S_2$ of codimension one. If $z_1\in B$, then we put $h=z_1$. Let $z_1\notin B$. If $z_3+\lambda z_1\in B$ with $\lambda\neq0$, then we put $h=\lambda^{-1}(z_3+\lambda z_1)$. Otherwise $z_3\in B$; take $b=z_4+\lambda z_1\in B$. Then $bz_3=-z_2+3\lambda z_3\in B$, so $z_2\in B$, and $z_2b=z_1-2\lambda z_2\in B$ gives $z_1\in B$; a contradiction. Thus $B$ contains an element $h$ with $L_h=L_{z_1}$, and $B$ is invariant under $L_{z_1}$. Since the eigenvalues of $L_{z_1}$ are pairwise distinct, $B$ is spanned by three of the elements $z_1,z_2,z_3,z_4$.

If $z_1\notin B$ (resp. $z_2\notin B$, $z_3\notin B$), then $z_2,z_4\in B$ (resp. $z_3,z_4\in B$; $z_2\in B$), and the equality $z_4z_2=-2z_1$ (resp. $z_4z_3=-z_2$; $z_2z_2=-3z_3$) gives a contradiction. Hence $B=\langle z_1,z_2,z_3\rangle$, and it follows from the multiplication table that this subspace is indeed a subalgebra. The theorem is proven.

\appendix
\section*{Appendix A. Differences from the published version}
\addcontentsline{toc}{section}{Appendix A. Differences from the published version}

\noindent\textbf{Theorem~\ref{theorem_left_com}:} the statement has been changed (the trivial annihilator is replaced by the trivial left annihilator; a counterexample to the original statement is given after the proof) and the proof corrected.

\noindent\textbf{Theorems~\ref{thm_codimW} and ~\ref{thm_codimW2}:} the proof has been corrected.

\noindent\textbf{Theorem~\ref{thm_codimS2}:} the statement has been changed (the subspace $\langle z_1,z_2,z_4\rangle$ is not a subalgebra, since $z_2z_2=-3z_3$) and the proof corrected.


\begin{thebibliography}{99}

\bibitem{Alison78} B.~Allison, \emph{A class of nonassociative algebras with involution containing the class of Jordan algebras}, Mathematische Annalen, 237 (1978), no.~2, 133--156; \doi{10.1007/BF01351677}.

\bibitem{AH81} B.~Allison and W.~Hein, \emph{Isotopes of some nonassociative algebras with involution}, Journal of Algebra, 69 (1981), no.~1, 120--142; \doi{10.1016/0021-8693(81)90132-0}.

\bibitem{BEC} E.~Barreiro, A.~Elduque, and C.~Mart\'inez, \emph{Derivations of the Cheng--Kac Jordan superalgebras}, Journal of Algebra, 338 (2011), no.~1, 144--156; \doi{10.1016/j.jalgebra.2011.02.016}; \arxiv{1101.0485}.

\bibitem{kac_can_10} N.~Cantarini and V.~Kac, \emph{Classification of linearly compact simple rigid superalgebras}, International Mathematics Research Notices, 2010 (2010), no.~17, 3341--3393; \doi{10.1093/imrn/rnp231}; \arxiv{0909.3100}.

\bibitem{dedkov89} A.~Dedkov, \emph{Some properties of quasiassociative and quasialternative algebras}, Siberian Mathematical Journal, 30 (1989), no.~3, 479--483; \doi{10.1007/BF00971505}.

\bibitem{dzhuma85} A.~Dzhumadil'daev, \emph{Simple Lie algebras with a subalgebra of codimension one}, Russian Mathematical Surveys, 40 (1985), no.~1, 215--216; \doi{10.1070/RM1985v040n01ABEH003537}.

\bibitem{Dzhuma01} A.~Dzhumadil'daev and S.~Abdykassymova, \emph{Leibniz algebras in characteristic $p$}, Comptes Rendus de l'Acad\'emie des Sciences. S\'erie I. Math\'ematique, 332 (2001), no.~12, 1047--1052; \doi{10.1016/S0764-4442(01)01941-3}.

\bibitem{Dzhuma07} A.~Dzhumadil'daev, \emph{Zinbiel algebras under $q$-commutators}, Journal of Mathematical Sciences (New York), 144 (2007), no.~2, 3909--3925; \doi{10.1007/s10958-007-0244-9}.

\bibitem{GR11} M.~Goze and E.~Remm, \emph{$2$-dimensional algebras}, African Journal of Mathematical Physics, 10 (2011), no.~1, 81--91.

\bibitem{kac77} V.~Kac, \emph{Classification of simple $\mathbb Z$-graded Lie superalgebras and simple Jordan superalgebras}, Communications in Algebra, 5 (1977), no.~13, 1375--1400; \doi{10.1080/00927877708822224}.

\bibitem{Kantor72} I.~Kantor, \emph{Certain generalizations of Jordan algebras}, Trudy Seminara po Vektornomu i Tenzornomu Analizu, 16 (1972), 407--499 (in Russian).

\bibitem{Kantor89term} I.~Kantor, \emph{An extension of a class of Jordan algebras}, Algebra and Logic, 28 (1989), no.~2, 117--121; \doi{10.1007/BF01979375}.

\bibitem{Kantor89tril} I.~Kantor, \emph{Terminal trilinear operations}, Algebra and Logic, 28 (1989), no.~1, 25--40; \doi{10.1007/BF01980606}.

\bibitem{Kantor89trudy} I.~Kantor, \emph{Some problems in $\mathcal{L}$-functor theory}, Trudy Instituta Matematiki (Novosibirsk), 16 (1989), 54--75 (in Russian).

\bibitem{Kantor90} I.~Kantor, \emph{The universal conservative algebra}, Siberian Mathematical Journal, 31 (1990), no.~3, 388--395; \doi{10.1007/BF00970345}.

\bibitem{kay} I.~Kaygorodov, \emph{$\delta$-Derivations of simple finite-dimensional Jordan superalgebras}, Algebra and Logic, 46 (2007), no.~5, 318--329; \doi{10.1007/s10469-007-0032-0}; (corrected version: \arxiv{1010.2419}).

\bibitem{kay_de2} I.~Kaygorodov, \emph{$\delta$-Derivations of classical Lie superalgebras}, Siberian Mathematical Journal, 50 (2009), no.~3, 434--449; \doi{10.1007/s11202-009-0049-9}; (corrected version: \arxiv{1010.2807}).

\bibitem{kay_lie2} I.~Kaygorodov, \emph{$\delta$-Superderivations of simple finite-dimensional Jordan and Lie superalgebras}, Algebra and Logic, 49 (2010), no.~2, 130--144; \doi{10.1007/s10469-010-9085-6}; (corrected version: \arxiv{1010.2423}).

\bibitem{kay_de5} I.~Kaygorodov, \emph{$\delta$-Superderivations of semisimple finite-dimensional Jordan superalgebras}, Mathematical Notes, 91 (2012), no.~1--2, 187--197; \doi{10.1134/S0001434612010208}; (corrected version: \arxiv{1106.2680}).

\bibitem{kay_okh} I.~Kaygorodov and E.~Okhapkina, \emph{$\delta$-Derivations of semisimple finite-dimensional structurable algebras}, Journal of Algebra and its Applications, 13 (2014), no.~4, 1350130; \doi{10.1142/S0219498813501302}; (corrected version: \arxiv{1212.0618}).

\bibitem{kaykantor} I.~Kaygorodov, \emph{On the Kantor product}, Journal of Algebra and its Applications, 16 (2017), no.~9, 1750167; \doi{10.1142/S0219498817501675}; (corrected version: \arxiv{1506.00736}).

\bibitem{kayvolk} I.~Kaygorodov and Yu.~Volkov, \emph{Conservative algebras of $2$-dimensional algebras, II}, Communications in Algebra, 45 (2017), no.~8, 3413--3421; \doi{10.1080/00927872.2016.1236935}; (corrected version: \arxiv{1507.02324}).

\bibitem{kuzmin89} E.~Kuzmin, \emph{The structure and representations of finite-dimensional Malcev algebras}, Quasigroups and Related Systems, 22 (2014), no.~1, 97--132.

\bibitem{shestakov_kuzm} E.~Kuzmin and I.~Shestakov, \emph{Non-associative structures}, in: A.~I.~Kostrikin and I.~R.~Shafarevich (eds.), Algebra VI, Encyclopaedia of Mathematical Sciences, 57, Springer, Berlin, 1995, 197--280.

\bibitem{Malcev55} A.~Malcev, \emph{Analytic loops}, Matematicheskii Sbornik (N.S.), 36(78) (1955), no.~3, 569--576 (in Russian).

\bibitem{GR2} H.~Petersson, \emph{The classification of two-dimensional nonassociative algebras}, Results in Mathematics, 37 (2000), no.~1--2, 120--154; \doi{10.1007/BF03322518}.

\bibitem{Popov2} A.~Popov, \emph{Differentiably simple Jordan algebras}, Siberian Mathematical Journal, 54 (2013), no.~4, 713--721; \doi{10.1134/S0037446613040113}.

\bibitem{shestakov2} A.~Shestakov, \emph{Ternary derivations of Jordan superalgebras}, Algebra and Logic, 53 (2014), no.~4, 323--348; \doi{10.1007/s10469-014-9293-6}.

\bibitem{wilansky75} A.~Wilansky, \emph{Subspaces, subalgebras and ideals of codimension one in complex algebras}, Journal of the London Mathematical Society (2), 9 (1974), no.~1, 87--92; \doi{10.1112/jlms/s2-9.1.87}.

\bibitem{zus10} P.~Zusmanovich, \emph{On $\delta$-derivations of Lie algebras and superalgebras}, Journal of Algebra, 324 (2010), no.~12, 3470--3486; \doi{10.1016/j.jalgebra.2010.09.032}; \arxiv{0907.2034}.

\end{thebibliography}
\end{document}